\documentclass[11pt,a4]{amsart}
\usepackage{amsthm,amssymb,amsmath,color,epsfig,comment}
\usepackage[colorlinks=true,allcolors=blue]{hyperref}

\newcommand{\ex}{{\rm e}\,}

\newcommand{\asy}{asymptotic}

\newcommand{\ts}{time series}

\renewcommand{\a}{\alpha}

\definecolor{darkblue}{rgb}{.1, 0.1,.8}
\definecolor{darkgreen}{rgb}{0,0.8,0.2}
\definecolor{darkred}{rgb}{.8, .1,.1}

\newtheorem{lemma}{Lemma}[section]

\newtheorem{theorem}[lemma]{Theorem}

\newcommand{\bbc}{{\mathbb C}}
\newtheorem{proposition}[lemma]{Proposition}
\newtheorem{definition}[lemma]{Definition}
\newtheorem{corollary}[lemma]{Corollary}
\newtheorem{example}[lemma]{Example}
\newtheorem{exercise}[lemma]{Exercise}
\newtheorem{remark}[lemma]{Remark}
\newtheorem{fig}[lemma]{Figure}
\newtheorem{tab}[lemma]{Table}

\newcommand{\bfQ}{{\bf Q}}

\newcommand{\bth}{\begin{theorem}}
\newcommand{\ethe}{\end{theorem}}

\newcommand{\bre}{\begin{remark}\em }
\newcommand{\ere}{\end{remark}}

\newcommand{\ble}{\begin{lemma}}
\newcommand{\ele}{\end{lemma}}
\newcommand{\sre}{stochastic recurrence equation}
\newcommand{\pp}{point process}
\newcommand{\bde}{\begin{definition}}
\newcommand{\ede}{\end{definition}}
\newcommand{\bco}{\begin{corollary}}
\newcommand{\eco}{\end{corollary}}

\newcommand{\bpr}{\begin{proposition}}
\newcommand{\epr}{\end{proposition}}

\newcommand{\bexer}{\begin{exercise}}
\newcommand{\eexer}{\end{exercise}}

\newcommand{\bexam}{\begin{example}}
\newcommand{\eexam}{\end{example}}

\newcommand{\bfi}{\begin{fig}}
\newcommand{\efi}{\end{fig}}

\newcommand{\btab}{\begin{tab}}
\newcommand{\etab}{\end{tab}}

\newcommand{\lhs}{left-hand side}
\newcommand{\fidi}{finite-dimensional distribution}
\newcommand{\rv}{random variable}

\newcommand{\sign}{{\rm sign}}

\newcommand{\var}{{\rm var}}

\newcommand{\cov}{{\rm cov}}

\newcommand{\as}{{\rm a.s.}}

\newcommand{\rhs}{right-hand side}
\newcommand{\df}{distribution function}

\newcommand{\dint}{\displaystyle\int}

\newcommand{\beao}{\begin{eqnarray*}}
\newcommand{\eeao}{\end{eqnarray*}\noindent}

\newcommand{\beam}{\begin{eqnarray}}
\newcommand{\eeam}{\end{eqnarray}\noindent}

\newcommand{\beqq}{\begin{equation}}
\newcommand{\eeqq}{\end{equation}\noindent}

\newcommand{\bce}{\begin{center}}
\newcommand{\ece}{\end{center}}

\newcommand{\barr}{\begin{array}}
\newcommand{\earr}{\end{array}}

\newcommand{\stp}{\stackrel{\P}{\rightarrow}}
\newcommand{\std}{\stackrel{d}{\rightarrow}}
\newcommand{\stas}{\stackrel{\rm a.s.}{\rightarrow}}

\newcommand{\stw}{\stackrel{w}{\rightarrow}}

\newcommand{\eqd}{\stackrel{d}{=}}

\newcommand{\vague}{\stackrel{\lower0.2ex\hbox{$\scriptscriptstyle
                    \it{v} $}}{\rightarrow}}
\newcommand{\weak}{\stackrel{\lower0.2ex\hbox{$\scriptscriptstyle
                    \it{w} $}}{\rightarrow}}
\newcommand{\what}{\stackrel{\lower0.2ex\hbox{$\scriptscriptstyle
                    \it{\hat{w}} $}}{\rightarrow}}

\newcommand{\bdis}{\begin{displaymath}}
\newcommand{\edis}{\end{displaymath}\noindent}

\renewcommand{\P}{\mathbb P}
\newcommand{\R}{\mathbb{R}}

\newcommand{\nto}{n\to\infty}
\newcommand{\kto}{k\to\infty}
\newcommand{\xto}{x\to\infty}

\newcommand{\ov}{\overline}
\newcommand{\wt}{\widetilde}
\newcommand{\wh}{\widehat}
\newcommand{\vep}{\varepsilon}

\newcommand{\regvary}{regularly varying}
\newcommand{\slvary}{slowly varying}
\newcommand{\regvar}{regular variation}

\newcommand{\bbr}{{\mathbb R}}

\newcommand{\bbz}{{\mathbb Z}}

\newcommand{\con}{convergence}

\newcommand{\evt}{extreme value theory}
\newcommand{\evd}{extreme value distribution}

\newcommand{\st}{such that}
\newcommand{\fif}{if and only if}
\newcommand{\wrt}{with respect to}

\newcommand{\fct}{function}

\newcommand{\ds}{distribution}

\newcommand{\rep}{representation}
\newcommand{\cmt}{continuous mapping theorem}
\newcommand{\seq}{sequence}

\newcommand{\pro}{probabilit}

\newcommand{\ms}{measure}

\newcommand{\ld}{large deviation}

\newcommand{\bfX}{{\bf X}}

\def\1{\ensuremath{\mathrm{1}\hspace{-.35em} \mathrm{1}}} % indicatrice

\def\E{{\mathbb E}}

\def\P{{\mathbb{P}}}
\def\R{\mathbb{R}}
\def\Z{\mathbb{Z}}
\renewcommand{\le}{\ensuremath{\leqslant}}
\renewcommand{\ge}{\ensuremath{\geqslant}}

\newcommand{\introo}[2]{{\left]{#1,\,#2\,}\right[\kern1pt}}

\newcommand{\intrfo}[2]{{\left[{#1,\,#2}\right[\kern1pt}}

\begin{document}
\today

\bibliographystyle{plain}
\title[Some variations on the extremal index]{Some variations on the extremal index}

\thanks{Thomas Mikosch's research is partially supported by Danmarks Frie Forskningsfond Grant No 9040-00086B}

\author[G. Buritica]{Gloria Buritic\'a}
\address{LPSM, Sorbonne Universit\'es\\
UPMC Universit\'e Paris 06, 
F-75005, Paris, France}
\email{gloria.buritica@sorbonne-universite.fr}

\author[N. Meyer]{Nicolas Meyer}
\address{Department  of Mathematics,
University of Copenhagen,
Universitetsparken 5,
DK-2100 Copenhagen,
Denmark}
\email{meyer@math.ku.dk}
\author[T. Mikosch]{Thomas Mikosch}
\address{Department  of Mathematics,
University of Copenhagen,
Universitetsparken 5,
DK-2100 Copenhagen,
Denmark}
\email{mikosch@math.ku.dk}
\author[O. Wintenberger]{Olivier Wintenberger}
\address{LPSM, Sorbonne Universit\'es\\
UPMC Universit\'e Paris 06, 
F-75005, Paris, France}
\email{olivier.wintenberger@upmc.fr}

\begin{abstract}{We re-consider Leadbetter's extremal index for 
stationary \seq s. It has interpretation as reciprocal of 
the expected size of an extremal cluster above high thresholds.
We focus on heavy-tailed \ts , in particular on \regvary\ stationary
\seq s, and discuss recent research in \evt\ for these models. 
A \regvary\ \ts\ has multivariate \regvary\ \fidi s. Thanks to results
by Basrak and Segers \cite{basrak:segers:2009} we have explicit 
\rep s of the limiting cluster structure of extremes, leading to explicit
expressions of the limiting \pp\ of exceedances and the extremal
index as a summary \ms\ of extremal clustering. 
The extremal index appears in various situations which do not seem
to be directly related, like the \con\ of maxima and  \pp es.
We consider different \rep s of the extremal index which arise from
the considered context.
We discuss the theory and 
apply it to a \regvary\ AR(1) process and the solution to an affine 
\sre .

}
\end{abstract}
\keywords{Extremal index, cluster Poisson process,
 extremal cluster, \regvary\ \ts , affine \sre ,
autoregressive process}
\subjclass{Primary 60G70 Secondary 60G55 60F99 60J10 62M10 62G32}
\maketitle

{\sc Some personal words by Thomas Mikosch.} During my PhD studies at
the University of Leningrad 1981--1984 I met Yasha Nikitin at 
seminars and workshops at LOMI (now POMI) and the university. I remember him 
as a professor who had a lot of humor and a balanced personality. 
Later, since the 1990s, Yasha Nikitin became a
representative of Russian Probability Theory and Mathematical Statistics 
with a high international reputation. I appreciated his cosmopolitan attitude. Whenever one needed
constructive advice as regards
some international scientific event 
(such as the European Meeting of Statisticians) 
or editorial issues, he would help. He supported 
the Bernoulli Society actively,  as an organization which embraces all European probabilists and statisticians. 
\par
I met Yasha at the Vilnius Conference in 2018 at the last time
and, as always,  I enjoyed his warm-hearted personality. He went from
us too early. His achievements for Probability Theory and Statistics
at the University of St. Petersburg and worldwide remain alive. 

\newpage

\section{Leadbetter's approach to modeling the extremes of a stationary \seq }\label{subsec:leadapproach}\setcounter{equation}{0}
The paper by Leadbetter \cite{leadbetter:1983} and  the 
book of Leadbetter, Lindgren and Rootz\'en \cite{leadbetter:lindgren:rootzen:1983}
provided a first systematic approach to the \evt\ of dependent stationary
\seq s. In particular, Leadbetter introduced mixing and anti-clustering
conditions, the conditions $D$ and $D'$,  which are tailored for the analysis of dependent extremal events. Moreover, \cite{leadbetter:lindgren:rootzen:1983} propagated the use of the
{\em extremal index} as a \ms\ for extremal clustering.   
\par
The idea of an extremal index originates from  
\cite{newell:1964,loynes:1965,obrien:1974}
who discovered that the maxima 
\beao
M_n=\max_{t=1,\ldots,n} X_t\,,\qquad n\ge 1\,,
\eeao
of numerous examples of dependent stationary \seq s $(X_t)$ with common \ds\ $F$ share the property 
that 
\beao
\P(M_n\le u_n)\approx \big[\P(X\le u_n)\big]^{n\,\theta_X}= \big((F(u_n))^n\big)^{\theta_X}\,,\qquad \nto\,,
\eeao 
for some number $\theta_X\in [0,1]$ provided $(u_n)$  is a \seq\ of high thresholds converging 
sufficiently fast to the 
right endpoint $x_F$ of $F$. Leadbetter
\cite{leadbetter:1983} made this notion precise 
as the {\em expected size of an extremal 
cluster of exceedances above high-level thresholds.} 
Since $(F(u_n))^n$ is the \df\ of the maximum of $n$ iid \rv s with common \ds\ $F$
at the threshold $u_n$, the quantity $\theta_X$ describes the shrinking effect
that the appearance of dependent extremes may have on the \ds\ of $M_n$
compared to $(F(u_n))^n$.
\par
Leadbetter defined the extremal index $\theta_X$ as follows:
assume that for every $\tau\in(0,\infty)$ there exists a \seq\ $(u_n(\tau))$
\st\ 
\beao
n\,\ov F(u_n(\tau))= n\,(1-F(u_n(\tau)))\to  \tau,
\eeao
 and 
there exists a number $\theta_X$ \st\
\beao 
\P(M_n\le u_n(\tau))\to
\ex^{-\tau\,\theta_X},\qquad  \nto\,.
\eeao
If such a number $\theta_X$ exists it belongs to the interval $[0,1]$ and is independent of the choice of the 
\seq s $(u_n)$.
\par
An immediate application is to the \con\ in \ds\ of the \seq\ 
$(M_n)$. Assume that
$(X_t)$ belongs to the maximum domain of attraction of an \evd\ $H$, i.e., for iid copies $(\wt X_t)$ 
of $X_1$, $\wt M_n=\max(\wt X_1,\ldots,\wt X_n)$, there exist constants $c_n>0$,
$d_n\in\bbr$\ \st\   $c_n^{-1}(\wt M_n-d_n)\std \xi$ as $\nto$ and 
$\xi$ has \ds\ $H$. Then if $(X_t)$ has an extremal index $\theta_X$ we have 
\beao
n\,\ov F(\underbrace{c_n\,x+d_n}_{=:u_n(\tau)})\to 
\underbrace{-\log H(x)}_{=:\tau}\,,\qquad \nto\,,\qquad x\in {\rm supp} H\,.
\eeao
and
\beao
\P\big(c_n^{-1}(M_n-d_n)\le x\big)\to H^{\theta_X}(x)\,,\qquad \nto\,,\qquad x\in 
{\rm supp} H\,. 
\eeao
\par
In the case of an iid \seq\ it is easily seen 
that $n\,\ov F(u_n(\tau))\to  \tau$ holds \fif\  $\P(M_n\le u_n(\tau))\to
\ex^{-\tau}$. Hence $\theta_X=1$. The extremal index $1$ is not exclusive to
iid \seq s. Indeed, in the book \cite{leadbetter:lindgren:rootzen:1983}
various examples of strictly stationary \seq s are considered for which
$\theta_X=1$. For example, if $(X_t)$ is a Gaussian stationary \seq\ 
whose autocovariance \fct\ satisfies $\cov(X_0,X_h)=o(1/\log h)$ as 
$h\to\infty$, then $\theta_X=1$.

\section{Sufficient conditions for the existence of the extremal index}

The extremal index is often interpreted as {\em the reciprocal of the 
expected size of an extremal cluster} for a stationary \seq\ $(X_t)$. 
We will give 
some justification for this interpretation.
\subsection{The method of block maxima}
The key is the definition of an 
{\em extremal cluster in the sample} $X_1,\ldots,X_n$: split
the sample into $k_n= [n/r_n]$ blocks of equal length $r_n$:
\beao
\underbrace{X_1,\ldots,X_{r_n}}_{\mbox{\small Block 1}}\,,\underbrace{X_{r_n+1},\ldots,X_{2\,r_n}}_{\mbox{\small Block 2}}\,,\ldots, \underbrace{X_{(k_n-1)\,r_n+1},\ldots,X_{k_n\,r_n}}_{\mbox{\small Block $k_n$}}\,,
\eeao
we ignore the last block of length less 
than $r_n$, and we
simply call a block an {\em extremal cluster} relative to
a high  threshold $u=u_n$ (this means that $u_n\uparrow x_F$ as $\nto$) 
if there is at least one exceedance of this threshold in this block. For
an \asy\ theory it will be important that $r=r_n\to\infty$  \st\ $r_n$ 
is small compared to $n$, i.e., $k_n\to \infty$. 
\par
In view of the stationarity of $(X_t)$ the {\em expected
cluster size of a block} is given~by
\beao
\E\Big[\sum_{t=1}^{r_n} \1(X_t>u_n)\,\Big|\, M_{r_n}>u_n\Big]&=&
\sum_{t=1}^{r_n} \dfrac{\P(X_t>u_n\,,M_{r_n}>u_n)}{\P(M_{r_n}>u_n)}\nonumber\\
&=&\sum_{t=1}^{r_n} \dfrac{\P(X_t>u_n)}{\P(M_{r_n}>u_n)}\nonumber\\&=&\dfrac{r_n\,\P(X>u_n)}{\P(M_{r_n}>u_n)}=:\dfrac1 {\theta_n}\,.
\eeao
Obviously, $\theta_n$ is a number in $[0,1]$.
Under mild regularity conditions the limit
$\theta=\lim _{\nto} \theta_n$ 
exists, assumes values in $(0,1]$ and coincides with Leadbetter's extremal index $\theta_X$; see Theorem~\ref{thm:1} below. For this reason, the extremal index $\theta_X$ is often referred to as {\em the reciprocal 
of the expected extremal cluster size above high thresholds.}

\begin{figure}[thbp]
\centerline{
\epsfig{figure=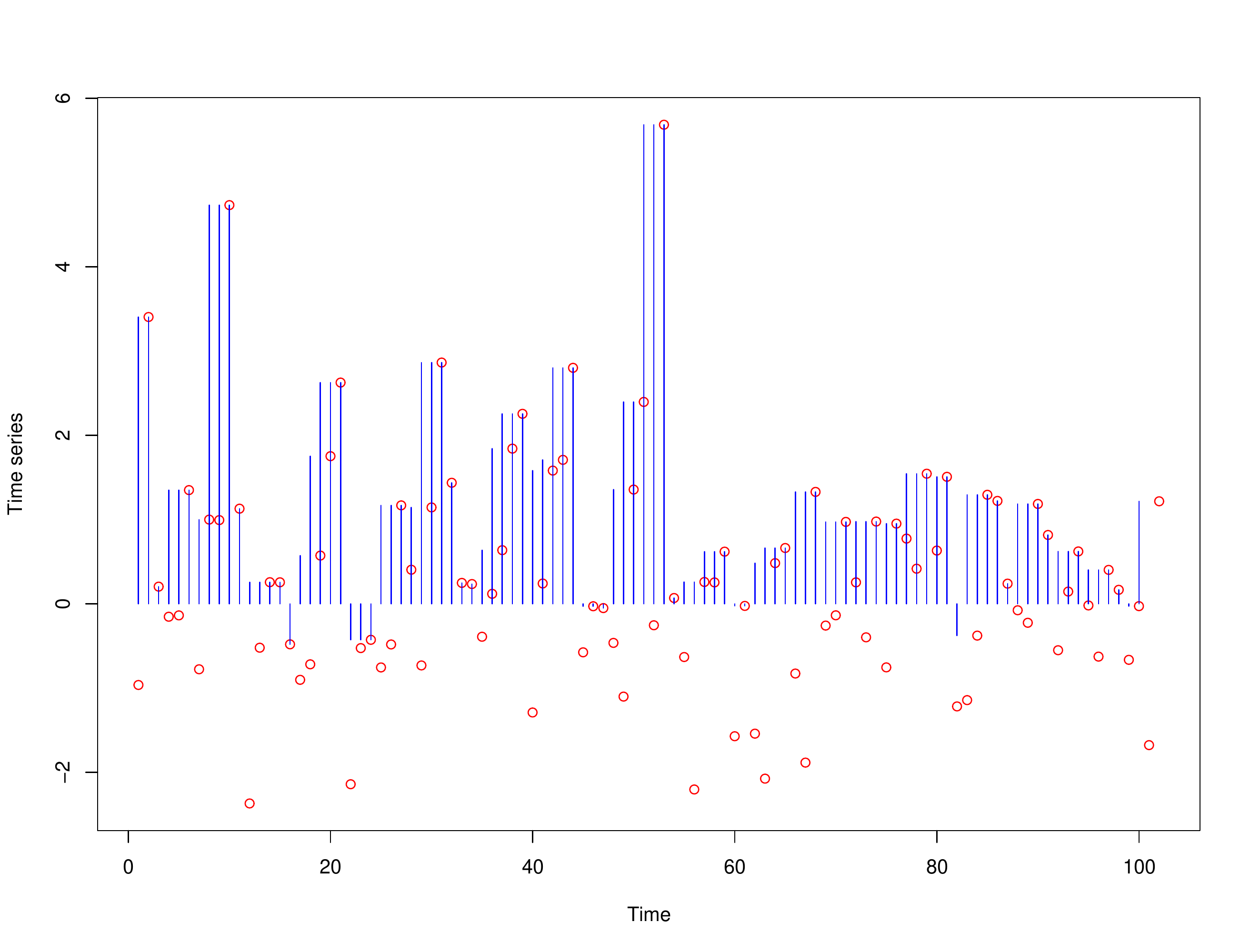,height=13cm,width=4cm,angle=-90}
}
\bfi{\small Visualization of the max-moving average $X_t=\max(Z_t,Z_{t+1},Z_{t+2})$,
$t=1,\ldots,100$, (blue)
for iid student noise $Z_t$, $t=1,\ldots,102$, with $\a=4$ degrees of freedom
(red dots). 
The values of $X_t$ typically appear in clusters of size 3. The process $(|X_t|)$
has extremal index $\theta_{|X|}=1/3$.
}\efi\label{fig:5c} 
\end{figure}
\begin{figure}[thbp]
\centerline{
\epsfig{figure=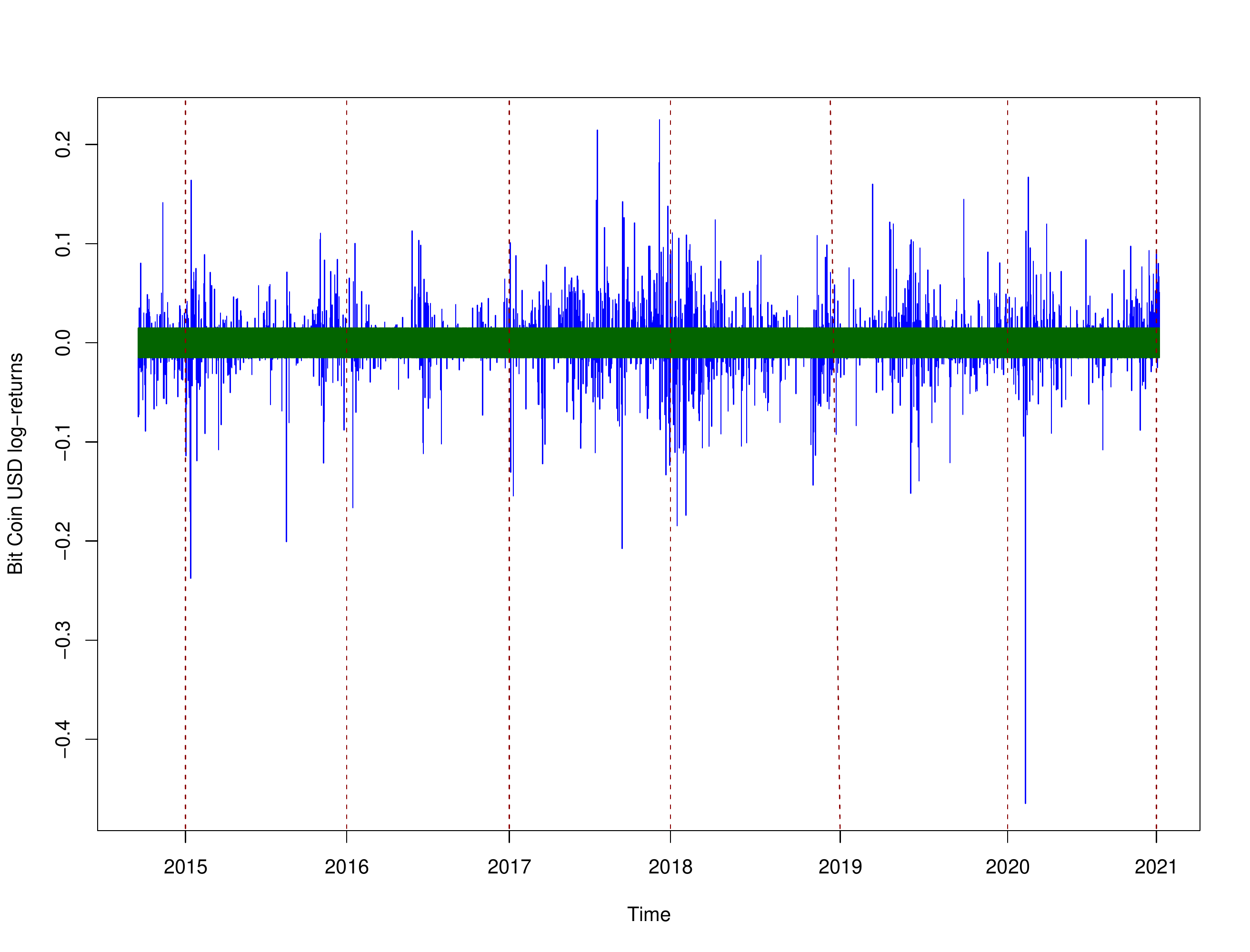,height=13cm,width=4cm,angle=-90}}
\bfi{\small The daily log-return series of the Bit Coin USD stock prices
from 17 September 2014 until 8 January 2021. We only show the returns
below -0.04 or above 0.04 which we interpret as extreme values.
These limits roughly correspond to the 10\% and 90\% quantiles of the data.
The extremes typically appear in clusters.}\label{fig:5a}\efi
\end{figure}

\subsection{Approximation of $\theta_X$ by $\theta_n$}
The following result can be found in slightly different forms in \cite{davis:hsing:1995}, proof of Lemma 2.8, \cite{segers:2005,basrak:segers:2009}.
\bth\label{thm:1}
Consider the following conditions:\\[2mm]
{\rm(1)} $(X_t)$ is a real-valued stationary \seq\ whose marginal \ds\ $F$
does not have an atom at the right endpoint $x_F$.\\
{\rm(2)} For a \seq\  $u_n\uparrow x_F$ and an
integer \seq\   $r=r_n\to\infty$ \st\ $k_n=[n/r_n]\to\infty$
the following {\em anti-clustering
  condition} is satisfied:
\beam\label{eq:71}
\lim_{k\to\infty}\limsup_{\nto} \P\big(M_{k,r_n}> u_n\,\mid\, X_0>u_n\big)=0\,.
\eeam
Here $M_{a,b}=\max_{i=a,\ldots,b}X_i$ for $a\le b$ such that $M_b = M_{a,b}$ with $a=1$.\\
{\rm(3)} A {\em mixing condition} holds: 
\beam\label{eq:mix}
\P(M_n\le u_n) - \big(\P(M_{r_n}\le u_n)\big)^{k_n}\to 0\,,\qquad \nto\,,
\eeam 
where $(u_n)$, $(k_n)$ and $(r_n)$ are as in $(2)$.\\
{\rm(4)} For all positive $\tau$ there exists a \seq\ $(u_n)=(u_n(\tau))$
\st\ $n\,\ov F(u_n)\to \tau$ and $(2),(3)$ are satisfied for these \seq s $(u_n)$.\\[2mm]
Then the following statements hold:\\[2mm]
{\rm 1.}
If $(1)$ and $(2)$ are satisfied  then
\beam\label{eq:73}
\lim_{k\to\infty}\limsup_{\nto}\big| 
\theta_n- \P\big(M_{k}\le u_n\,\mid\, X_0>u_n\big)\big| =0\,,
\eeam
and $\liminf_{\nto} \theta_n>0$.\\
{\rm 2.}
If $(1)$ and $(4)$ are satisfied and $\theta=\lim_{\nto}\theta_n$ exists, 
then $\theta_X\in (0,1]$ exists and coincides with $\theta$.
\ethe
Condition \eqref{eq:mix} is satisfied for strongly mixing $(X_t)$
with mixing rate $(\a_h)$ if one can find integer \seq s $(\ell_n)$ and 
$(r_n)$ \st\ $\ell_n/r_n\to 0$, $r_n/n\to 0$ and $k_n\a_{\ell_n}\to0$ 
as $\nto$.
Anti-clustering conditions are common in \evt\
since Leadbetter introduced the $D'$ condition which is much stronger than \eqref{eq:71} but also easily verified on examples. 
The goal of such a condition is to avoid that the stationary \seq\ stays above a high threshold for too long. 
\par
Relation \eqref{eq:73} is in agreement with O 'Brien's \cite{obrien:1987}
characterization of the extremal index of $(X_t)$ as the limit 
\beam\label{eq:obrien}
\theta_X=\lim_{\nto}\P(M_{\ell_n}\le u_n\,\mid \,X_0>u_n),
\eeam 
for a \seq\ $(\ell_n)$ with $\ell_n/n\to 0$, thresholds $u_n\uparrow x_F$
\st\ $n\,\ov F(u_n)\to 1$ as $\nto$, provided a mixing condition holds. O'Brien's condition \eqref{eq:obrien} has the advantage that it avoids
the definition of an extremal cluster. 
\par
\bre\label{rem:1}
Relation \eqref{eq:73} provides a constructive way of calculating $\theta_X$:
if we know that the limits $f(k):=\lim_{\nto}\P\big(M_{k}\le u_n\,\mid\, X_0>u_n\big)$ exist for every $k\ge 1$ then we can try to 
derive $\theta_X=\lim_{\kto}f(k)$. In Section~\ref{sec:2} we will follow 
this approach in the case of a {\em \regvary\ \seq }.
\ere

\section{Regularly varying \seq s}\label{sec:2}\setcounter{equation}{0}
\subsection{Definition and examples}
As a matter of fact, clusters of extremes are more prominent in 
stationary \seq s with heavy-tailed marginal \ds . To illustrate this fact,
consider a stationary causal AR(1) process which solves the difference equation 
$X_t=\varphi\,X_{t-1 }+Z_t$, $t\in\bbz$, for an iid noise \seq\ $(Z_t)$.
Necessarily, $\varphi\in (-1,1)$ and, if $(Z_t)$ is iid standard normal then
$(|X_t|)$ has extremal index $\theta_{|X|}=1$ (see \cite{leadbetter:lindgren:rootzen:1983}), while for iid student noise 
$(Z_t)$ with $\alpha$ degrees of freedom we have $\theta_{|X|}=1-|\varphi|^\a$;
see Example~\ref{exam:may1a} below.
Thus, the smaller $\a$ (the heavier the tail) for given $\varphi$ the closer $\theta_{|X|}$ to zero.
\par
An AR(1) process with student noise is an example of a {\em \regvary\ \ts }. This class
of heavy-tailed processes has been studied rather extensively in the last 15 years; 
see \cite{resnick:2007} for some basics about multivariate \regvar , and \cite{kulik:soulier:2020} for a recent textbook treatment.
This class was considered in full generality first by 
\cite{davis:hsing:1995}: they required that the \fidi s of the 
process satisfy a multivariate \regvar\ condition; see 
\cite{resnick:1987,resnick:2007} for the definition of this notion.
It is an extension of power-law tail behavior from the univariate to the 
multivariate case defined via the vague \con\ of tail \ms s
with infinite limit \ms s which have the homogeneity property. 
\par
Here we will follow an alternative approach by \cite{basrak:segers:2009} tailored for stationary \seq s, 
avoiding the vague \con\ concept. They proved 
that a real-valued stationary \seq\ $(X_t)$ is \regvary\ with index $\a>0$
in the sense of \cite{davis:hsing:1995}
\fif\ there exists a \seq\ $(\Theta_t)$ and a Pareto$(\a)$ distributed 
\rv\ $Y$, i.e., $\P(Y>y)=y^{-\a}$, $y>1$, \st\ $(\Theta_t)$ and $Y$
are independent and,
for all $h\ge 0$,
\beao
\P\big(x^{-1}(X_t)_{|t|\le h }\in \cdot\,\big|\,|X_0|>x\big)\stw
\P\big(Y\,(\Theta_t)_{|t|\le h}\in \cdot\big)\,,\qquad \xto\,.
\eeao
In the latter relation $x$ can be replaced by any \seq\ $(a_n)$
\st\ $n\,\P(|X|>a_n)\to 1$ 
as $\nto$. 
Moreover, by definition, $|\Theta_0|=1$ a.s. 
The \seq\ $(\Theta_t)$ is the {\em spectral tail process} of the \regvary\
process $(X_t)$; it describes the propagation  of a value 
$|X_0|>x$ for large $x$ through the stationary \seq\ $(X_t)$ into its past and future.
\bexam\label{exam:ar1} \rm We consider a stationary AR(1) process given as the causal solution 
to the difference equation $X_t=\varphi\,X_{t-1}+Z_t$, $t\in\bbz$,
where $(Z_t)$ is iid \regvary\ with index $\a$ (e.g. Pareto$(\a)$ or student$(\a)$). This means that a generic element $Z$ satisfies $\lim_{\xto}\P(\pm Z>x)/\P(|Z|>x)=p_\pm$ for non-negative values $p_\pm$ \st\ $p_++p_-=1$, and $\P(|Z|>x)=L(x)x^{-\a}$, $x>0$,
for some \slvary\ \fct\ $L$. Then a generic element $X$ inherits the \regvary\
tail behavior from $Z$ (see \cite{davis:resnick:1985}):
\beao
\dfrac{\P(\pm X>x)}{\P(|Z|>x)}\sim \sum_{j=0}^\infty\big[ p_\pm\,(\varphi^j)_\pm ^{\a}
+p_\mp \,(\varphi^j)_{\mp}^{\a} \big] = \P(\Theta_0=\pm 1)(1-|\varphi|^\alpha)\,.
\eeao
But even more is true: $(X_t)$ is a \regvary\ \ts\
with spectral tail process
\beam\label{eq:march11a}
\Theta_t=\Theta_Z\,\sign(\varphi^{J+t})\,|\varphi|^{t} \1(J+t\ge 0)=\Theta_0\,\varphi^t\,\1(J+t\ge 0)\,,\qquad t\in\bbz\,,\nonumber\\ 
\eeam
where
$\P(\Theta_Z=\pm 1)=p_\pm$, $\Theta_Z$ is independent of $J$ which has \ds
\beao
\P(J=j)= (1-|\varphi|^\a)\,|\varphi|^{j\,\a}\,,\qquad j\ge 0\,.
\eeao
In particular, the forward spectral tail process is given by $\Theta_t=\Theta_0\,\varphi^t$, $t\ge 0$.
\eexam
\bexam\label{exam:2}\rm We consider the unique causal solution to the affine \sre\
$X_t= A_t\,X_{t-1}+B_t$, $t\in\bbz$, for an iid \seq\ $((A_t,B_t))_{t\in\bbz}$
with generic element $(A,B)\in \bbr_+^2$. We assume that the \ds\ of $(A,B)$
satisfies the conditions of the Kesten-Goldie theory; see 
\cite{kesten:1973,goldie:1991}, cf. \cite{buraczewski:damek:mikosch:2016} for a textbook treatment.
The most important condition in this context is the existence
of a unique solution $\a>0$ to the equation $\E[A^\a]=1$ which yields the 
tail index $\a$. Under these conditions for a generic element $X$, there exists a positive constant $c_+$ such that
\beao
\P( X> x)\sim c_+ \,x^{-\a}\,,\qquad \xto\,.
\eeao
The forward spectral process is then
given by 
\beao
(\Theta_t)_{t\ge 0}= (\Pi_t)_{t\ge 0}\,, \mbox{ where $\Pi_t=A_1\cdots A_t$}\,,
\eeao
while the backward spectral tail process $(\Theta_t)_{t\le -1}$ has a rather 
complicated structure.
\par
Writing $S_t=\log \Pi_t=\sum_{i=1}^t \log A_i$, $t\ge 1$, we observe that  
$(S_t)$ constitutes  a random walk with a negative drift. Indeed, by Jensen's 
inequality we have $\E[\log (A^\a)]< \log (\E[A^\a])=0$. 
\eexam
\subsection{The extremal index}
Following Remark~\ref{rem:1}, we will derive the extremal index $\theta_X$
of a stationary non-negative \regvary\ \seq\ $(X_t)$ in terms of its spectral tail process.
First, we observe that by virtue of the \cmt , as $\nto$ for $k\ge 1$, 
\beao\lefteqn{
\P\big(a_n^{-1} M_k\le 1\,\big|\,X_0>a_n\big)}\\&\to &
\P\big(Y\,\max_{1\le t\le k}\Theta_t \le 1\big)
=\P\big(\max_{1\le t\le k}\Theta_t^\a \le Y^{-\a}\big)\\
&=& \E\big[\big(1- \max_{1 \le t\le k} \Theta_t^\a\big)_+\big]
= \E\big[\max_{0\le t\le k} \Theta_t^\a- \max_{1\le t\le k} \Theta_t^\a\big]\,.
\eeao 
Here we used the fact that $Y^{-\a}$ is $U(0,1)$ uniformly distributed and $\Theta_0=1$ a.s. Using dominated \con\ as $\kto$, we proved under the 
anti-clustering condition \eqref{eq:71} that
\beao
\lim_{\nto}\theta_n&=&\lim_{\kto}\lim_{\nto} \P\big(a_n^{-1} M_k\le 1\big|\,X_0>a_n\big)\\
&=& \lim_{\kto} 
  \E\big[\max_{0\le t\le k} \Theta_t^\a- \max_{1\le t\le k} \Theta_t^\a\big]\\
&=& \E\big[\big(1- \max_{t\ge 1} \Theta_t^\a\big)_+\big]\,.
\eeao
From Theorem~\ref{thm:1}
we obtain the following result in \cite{basrak:segers:2009}.
\bco\label{cor:bstheta} 
Consider a non-negative stationary \regvary\ process 
$(X_t)$ with index $\alpha>0$. Then the following statements hold:\\[2mm]
{\rm 1.}
If the anti-clustering condition 
\eqref{eq:71} holds for $(u_n)=(x\,a_n)$ and some $x>0$ then the limit
$\theta=\lim_{\nto} \theta_n$ exists, is positive and  has the \rep s
\beam\label{eq:march12a}
\theta&=& \P\big(Y\,\sup_{t\ge 1} \Theta_t\le 1\big)
= \E\big[\big(1-\sup_{t\ge 1} \Theta_t^\alpha\big)_+\big]
=\E\big[\sup_{t\ge 0} \Theta_t^\alpha -\sup_{t\ge 1} \Theta_t^\alpha\big]\,.\nonumber\\
\eeam
{\rm 2.} 
If \eqref{eq:71} and the mixing condition 
\eqref{eq:mix} hold for $(u_n)=(x\,a_n)$ 
and all $x>0$  then the extremal index $\theta_X$ exists and 
coincides with~$\theta$.
\eco
The \rep s of $\theta$ given in \eqref{eq:march12a} only depend 
on the forward spectral process $(\Theta_t)_{t\ge 0}$. In Proposition~\ref{cor:cluster1} below
we provide \rep s of the extremal index $\theta_{|X|}$ depending
on the whole spectral tail process $(\Theta_t)_{t\in\bbz}$.

\bexam\label{exam:may1a} \rm We consider the \regvary\ AR(1) process from Example~\ref{exam:ar1}. 
It can be shown to satisfy the anti-clustering and mixing conditions of
Theorem~\ref{thm:1}. We conclude from Corollary~\ref{cor:bstheta} and the form of the 
spectral tail process given in \eqref{eq:march11a}
that
\beao
\theta_{|X|}=\E\big[\big(1- \max_{t\ge 1} \Theta_t^\a\big)_+\big]=
1- \max_{t\ge 1} |\varphi|^{\a\,t}=1-|\varphi|^\a\,. 
\eeao
This formula was already achieved in \cite{davis:resnick:1985}
in the wider context of linear processes.
\eexam
\bexam\rm  We consider the \regvary\ solution of an affine \sre\ under the 
conditions and with the notation of Example~\ref{exam:2}. It can be shown to satisfy the anti-clustering and mixing conditions of
Theorem~\ref{thm:1}; see \cite{buraczewski:damek:mikosch:2016}.  We conclude from this result that $(X_t)$ has extremal index
\beao
\theta_{X}= 
\E\big[\big(1- \max_{t\ge 1} \Pi_t \big)_+\big]=
\E\big[\big(1- \exp\big( \max_{t\ge 1} S_t\big) \big)_+\big]\,,
\eeao
where $S_t=\sum_{i=1}^t \log A_i$, $t\ge 1$, is a random walk with a negative drift. This value of $\theta_X$
was derived in \cite{haan:resnick:rootzen:vries:1989}.
In that paper a Monte Carlo simulation procedure for the evaluation of  
$\theta_X$ was proposed.
Direct calculation of $\theta_X$ is difficult; see Example \ref{ex:br} for an exception. \eexam
\subsection{The extremal index and \pp\ \con\ toward a cluster Poisson process}
\subsubsection{A useful auxiliary result}
\ble\label{lem:1}
Consider a non-negative stationary \regvary\ \seq\ $(X_t)$ with index $\a>0$ and assume that
\eqref{eq:71} holds for $(u_n)=(x\,a_n)$ and all $x>0$.
Then 
\beao
\|\Theta\|_\a^\a:=\sum_{j\in\Z}\Theta_j^\alpha<\infty\qquad \as \,
\eeao
In particular, $\Theta_t\to 0$ a.s. as $|t|\to\infty$, and the 
time $T^\ast$ of the largest record  
 of $(\Theta_t)$ is finite, i.e., $|T^\ast|$ is the smallest integer \st
\beao
\Theta_{T^\ast}= \max_{t\in \bbz} \Theta_t\,. 
\eeao
\ele
\begin{proof} Write $(Y_t)= Y\,(\Theta_t)$ where the Pareto$(\a)$ variable $Y$ and the spectral tail process $(\Theta_t)$ are independent.
 We start by showing 
\beam\label{eq:june25}
Y_t\stas  0\,,\qquad t\to\infty\,.
\eeam
Since $(X_t)$ is \regvary\  we have for all $x>0$ and integers $k\ge 1$,
\beao
\lim_{h\to\infty}\lim_{\nto}\P\big(M_{ k,k+h}>x\,a_n\;\mid\; X_0>a_n\big)&=& 
\lim_{h\to\infty}\P\big(\max_{  k\le t\le k+h} Y_t>x\big)\\&=& \P\big(\max_{  t\ge k} Y_t>x\big)\,.
\eeao
On the other hand,  
using the anti-clustering condition
\eqref{eq:71} for all $x\in (0,1]$, we have for fixed $k,h\ge 1$,
\beao\lefteqn{
\lim_{\nto}
\P\big(M_{  k,k+h}>x\,a_n\;\mid\; X_0>a_n\big)}\\
&\le &\limsup_{\nto}\P\big(M_{k,r_n}>x\,a_n\;\mid\;X_0>x\,a_n\big)\dfrac{\P(X>x\,a_n)}{\P(X>a_n)}\\
&=&x^{-\a}\,\limsup_{\nto}\P\big(M_{k,r_n}>x\,a_n\;\mid\;X_0>x\,a_n\big)
\\&=&x^{-\a}\varepsilon_k\,,
\eeao
and the \rhs\ term $\varepsilon_k$ vanishes for large $k$.
Hence, letting $h\to \infty$, we obtain for all $x>0$,
\beao
\P\big(\max_{  t\ge k} Y_t>x\big)\le x^{-\a}\varepsilon_k\,,
\eeao
 and therefore
\beao
\lim_{\kto} \P\big(\max_{t\ge k} Y_t>x\big)\le \lim_{\kto} x^{-\a}\varepsilon_k
=0\,,
\eeao
implying $\max_{t\ge k} Y_t\stp 0$ as $\kto$.
Since $(Y_t)= Y\,(\Theta_t)$ a.s. and $ Y>0$ is independent of $(\Theta_t)$ this
is only possible if $\max_{t\ge k} \Theta_t\stp 0$ as $\kto$ but the latter relation
is equivalent to  $\Theta_t\stas 0$ as $t\to\infty$, implying \eqref{eq:june25}.\\[2mm]
Next we show that 
\beao
Y_{-t}\stas 0\,,\qquad t\to\infty\,.
\eeao
Since $Y_t\stas 0$  as $t\to\infty$ and $Y_0>1$ a.s. the following relation
holds 
\beao
%\P\Big(\bigcup_{i\ge 0} \Big\{Y_i\ge \vep > \max_{t> i}Y_t\Big\}\Big)=\sum_{i\in \Z}\P\Big(Y_i\ge \vep > \max_{t> i}Y_t\Big)&\le&1\,,\qquad \vep>0\, ,\\
\P\big(\bigcup_{i\ge 0} \big\{Y_i\ge 1 > \max_{t> i}Y_t\big\}\big)=
\sum_{i\ge 0}\P\big(Y_i\ge1 > \max_{t> i}Y_t\big)=1\,.
\eeao 
Suppose that
$\P(\sum_{j\le 0}\1(Y_j>\vep)=\infty)>0$ for some $\vep>0$. Then there exists some $i\ge 0$ \st\
\beao
\P\Big(\sum_{j\le 0}\1(Y_j>\vep)=\infty,Y_i\ge 1 > \max_{t> i}Y_t\Big)>0\,. 
\eeao
%If such an $i$ did not exist
%we would have 
%\beao
%0&=&\sum_{i\ge 0} \P\Big(\sum_{j\le 0}\1(Y_j>\vep)=\infty, Y_i\ge 1 > \max_{t> i}Y_t\Big)\\
%&=&\P\Big(\sum_{j\le 0}\1(Y_j>\vep)=\infty\Big)\,,
%\eeao
%yielding a contradiction.
We recall the time-change formula from \cite{basrak:segers:2009}:
\beam\label{eq:tcp}
\P((\Theta_{-h},\ldots,\Theta_{h})\in \cdot \mid \Theta_{-t}\ne {0})&=&\E\big[\dfrac{\Theta_{t}^\alpha}{\E\big[\Theta_{t}^\alpha\big]}
\1\big(\dfrac{(\Theta_{t-h},\ldots,\Theta_{t+h})}{\Theta_{t}}\in \cdot\big)\big]\,.\nonumber\\
\eeam
In particular, $\P(\Theta_t\ne 0)=\E[\Theta_{t}^\alpha]=1$ \fif\ for all $h\ge 0$,
\beao
\P((\Theta_{-h},\ldots,\Theta_{h})\in \cdot  )=\E\big[
\dfrac{\Theta_{t}^\alpha}{\E\big[\Theta_{t}^\alpha\big]}\,\1
\big(\dfrac{(\Theta_{t-h},\ldots,\Theta_{t+h})}{\Theta_{t}}\in \cdot\big)\big]\,.
\eeao
Therefore
\beao 
\infty &=& \E\big[\sum_{j\le 0} \1(Y_j>\vep)\,\1\big(Y_i\ge 1 > \max_{t> i}Y_t\big)\big]\\
&=&\sum_{j\le 0} \P \big(Y_j>\vep, Y_i\ge 1 > \max_{t> i}Y_t\big)\\
&=&\sum_{j\le 0} \int_1^\infty\E\big[\1\big(y \,\Theta_j>\vep\,,y \,\Theta_i\ge 1 > y\, \max_{t> i}\Theta_t\big)\big]d\big(-y^{-\a}\big)\\
&=&\sum_{j\le 0} \int_1^\infty\E\big[\Theta_{-j}^\a\,\1\big(y >\vep\,\Theta_{-j}\,, 
y \,\dfrac{\Theta_{i-j}}{\Theta_{-j}}\ge 1 > y \max_{t> i-j}\dfrac{\Theta_t}{\Theta_{-j}}\big)\big]d\big(-y^{-\a}\big)\\
%&=&\sum_{j\le 0} \E\Big[\int_0^\infty\Theta_{-j}^\a\1\Big( y >\vep\,\Theta_{-j},y \,\Theta_{i-j}\ge \Theta_{-j} > y \max_{t> i-j}\Theta_t\Big)d\big(-y^{-\a}\big)\Big]\\
&\le& \vep^{-\a}\sum_{j\le 0} \E\big[\int_1^\infty\1\big(z>1\,,z \,\Theta_{i-j}\ge \vep^{-1} > z \max_{t> i-j}\Theta_t\big)d\big(-z^{-\a}\big)\big]\\
&=& \vep^{-\a} \sum_{j\le 0}\P\big(Y_{i-j}\ge \vep^{-1} >  \max_{t> i-j} Y_t\big) \\
&=&\vep^{-\a}\sum_{k\ge i} \P\big(Y_k\ge  \vep^{-1} > \max_{t> k}Y_t\big)\\
&\le&\vep^{-\a}\,. 
\eeao
In the last step we used the fact that the events  $\{Y_k\ge \vep^{-1} > \max_{t> k}Y_t\}$, $k\ge i$, are disjoint.
Thus we got a contradiction. 
This proves that for all $\vep>0$ there exist only finitely many 
$j\le 0$ such that $Y_j>\vep$, hence $Y_t\stas 0$ 
and also $\Theta_t\stas 0$  as $t\to -\infty$,
as desired.
\par
In particular, the time   $T^\ast$ of the largest record  of the \seq\ 
$(\Theta_t)$ is finite a.s.
\par
Now suppose that
$\P(\sum_{j\in \Z}\Theta_j^\alpha=\infty)>0$. 
Then there exists an $i\in \Z$ such that 
\beao
\P\big(\sum_{j\in \Z}\Theta_j^\alpha=\infty\,, T^\ast=i\big)>0\,,
\eeao
and an application of the time-change formula \eqref{eq:tcp} yields
\beao
\infty &= &\E\Big[\sum_{j\in \Z}\Theta_j^\alpha \1(T^\ast=i)\Big]= \sum_{j\in \Z}\E\big[\Theta_j^\alpha \1(T^\ast=i)\big]\\
&= &\sum_{j\in \Z} \P(T^\ast=i-j)
=1\,,
\eeao
leading to a contradiction. Thus $\sum_{j\in \Z}\Theta_j^\alpha<\infty$ a.s.
This proves the lemma. 
\end{proof}
\subsubsection{Point process \con\ toward cluster Poisson processes}
The following \pp\ result was proved in \cite{davis:hsing:1995}
and re-proved in \cite{basrak:segers:2009} by using the 
terminology of the spectral tail process.
\par
We adapt the mixing condition in \cite{davis:hsing:1995} 
tailored for \pp\ \con . It is expressed in terms of the 
Laplace \fct als of \pp es. Recall that a \pp\ $N$ with state space
$E=\bbr_0=\bbr\backslash \{0\}$
has Laplace \fct al 
\beao
\Psi_N(g)=\E\left[\exp\left(-\dint_E g\,dN\right)\right]\qquad \mbox{ for  
$g\in\bbc_K^+$,}
\eeao 
where the set $\bbc_K^+$ consists of the continuous \fct s on $E$ with 
compact support. Since $0$ is excluded from $E$ this means that $g\in \bbc_K^+$
 vanishes in some neighborhood of the origin.  
Moreover, we have the weak \con\ of \pp es $N_n\std N$ on $E$ \fif\
$\Psi_{N_n}\to \Psi_N$ pointwise; see \cite{resnick:1987,resnick:2007}.\\[1mm]
{\bf Mixing condition ${\mathcal A}(a_n)$}
Consider 
integer \seq s $r_n\to\infty$ and $k_n=[n/r_n]\to\infty$
and the \pp es with state space $E=\bbr_0$,
\beao
N_n=\sum_{i=1}^n \vep_{a_n^{-1}X_i}\qquad\mbox{ and } 
\qquad \wt N_{r_n}=\sum_{i=1}^{r_n} \vep_{a_n^{-1}X_i}\,,\qquad  n\ge 1\,,
\eeao
where $\vep_x$ denotes Dirac \ms\ at $x$. 
The stationary \regvary\ \seq\ $(X_t)$ satisfies ${\mathcal A}(a_n)$ if 
there exist $(r_n)$ and $(k_n)$
\st
\beam\label{eq:tod}
\Psi_{N_n}(g) - \big(\Psi_{\wt N_{r_n}}(g)\big)^{k_n}\to 0\,,\qquad\nto\,,\quad  g\in\bbc_K^+\,.
\eeam
\bre
This condition is satisfied for a strongly mixing \seq\ $(X_t)$
with mixing rate $(\a_h)$ if one can find integer \seq s $(\ell_n)$ and 
$(r_n)$ \st\ $\ell_n/r_n\to 0$, $r_n/n\to 0$ and $k_n\a_{\ell_n}\to0$. This is a 
very mild condition indeed.
Relation~\eqref{eq:tod} 
ensures that, if $N_n\std N$ on the state space $E$, then also  
$\sum_{i=1}^{k_n}\wt N_{r_n}^{(i)}\std N$ 
where $(\wt N_{r_n}^{(i)})_{i=1,\ldots,k_n}$ are iid
copies of $\wt N_{r_n}$. This fact ensures that the limit processes considered
are infinitely divisible; cf. \cite{kallenberg:2017}.
\ere

\bth\label{thm:davishsing2}
Consider a stationary \regvary\ sequence $(X_t)$ with index
$\alpha>0$. We assume the following conditions:\\[1mm]
{\rm (1)} The mixing condition ${\mathcal A}(a_n)$ for integer \seq s
$r_n\to\infty$ \st\ $k_n=[n/r_n]\to\infty$ as $\nto$.\\
{\rm (2)} The anti-clustering condition \eqref{eq:mix}
for the same \seq\ $(r_n)$\,.\\[1mm]
Then we have the \pp\ \con\
 on the state space $\bbr_{0}$ 
\beam\label{eq:repcpp}
N_n= \sum_{i=1}^n\vep _{a_n^{-1} X_i}\std N=\sum_{i=1}^\infty\sum_{j=-\infty}^\infty\vep_{\Gamma_i^{ -1/\a} Q_{ij}}\,,
\eeam
where\\[1mm]
$\bullet$~$\sum_{j=-\infty}^\infty\vep_{ Q_{ij}}$, $i=1,2,\ldots$,  
is an iid \seq\ of \pp es with state space $\bbr$.  A generic element 
$Q=(Q_j)$ of the  \seq\  $Q^{(i)}=(Q_{ij})_{j\in\bbz}$,  $i=1,2,\ldots$, 
has the \ds\ of the {\em spectral cluster process} 
\beao
Q=\Big(\dfrac{\Theta_t}{\|\Theta\|_\alpha}\Big)_{t\in\bbz}\,.
\eeao
$\bullet$~$(\Gamma_i)$ are the points of a unit rate 
homogeneous Poisson process on $(0,\infty)$. 
$\bullet$~$(\Gamma_i)$ and $(Q^{(i)})_{i=1,2,\ldots}$ are independent.
\ethe
\bre
In view of Lemma~\ref{lem:1} we know that 
$\|\Theta\|_\a<\infty$ a.s. Hence the spectral cluster 
process $Q$ is well defined.
\ere
Since the Poisson points $(\Gamma_i^{-1/\a})$ and the \seq\ of iid \pp es
$\big(\sum_{j\in\bbz}\vep_{Q_{ij}}\big)$ are independent it is not difficult
to calculate the Laplace \fct al of the limit process $N$:
\beao
\Psi_N(g)= \exp\Big(-\int_0^\infty \E\big[1- \ex^{-\sum_{j\in\bbz} g(y\,Q_j)}  
\big]\,d(-y^{-\a})\Big)\,,\qquad g\in\bbc_K^+\,.
\eeao
%where $Q=(Q_j)_{j\in\bbz}$ is a generic element of the iid \seq\ $(Q^{(i)})$.
Now we apply the change of variables
$z= y\,|Q_{T^\ast}|$ in $\Psi_N(g)$ where 
\beao
|Q_{T^\ast}|= \dfrac{|\Theta_{T^\ast}|}{\|\Theta\|_\a}
= \dfrac{\max_{t\in\bbz}|\Theta_t|}{\Big(\sum_{j\in\bbz} |\Theta_j|^\a\Big)^{1/\a}}\,.
\eeao 
Then we obtain for $g\in\bbc_K^+$,
\beao
\Psi_N(g)&=& \exp\Big( -\E[|Q_{T^\ast}|^\a]\\&&\times
\int_0^\infty 
\E\Big[
\dfrac{|Q_{T^\ast}|^\a}{\E[|Q_{T^\ast}|^\a]}
\big(1- 
\ex^{-\sum_{j\in\bbz} g(z\,Q_j/|Q_{T^\ast}|)} \big) 
\Big]\,d(-z^{-\a})
\Big)\,.
\eeao
According to Proposition~\ref{cor:cluster1} below, 
$\theta_{|X|}=\E[|Q_{T^\ast}|^\a]$. Now, changing the \ms\ 
with the density $|Q_{T^\ast}|^\a/\E[|Q_{T^\ast}|^\a]$ and 
writing $\wt Q=(\wt Q_j)_{j\in\bbz}$ for the \seq\ $Q/|Q_{T^\ast}|$
under the new \ms , 
we arrive at
\beao
\Psi_N(g)&=& \exp\Big( -
\int_0^\infty 
\E\big[\big(1- 
\ex^{-\sum_{j\in\bbz} g(z\,\wt Q_j|)} \big) 
\big]\,d\big(-(z/\theta_{|X|}^{1/\a}\big)^{-\a}\big)
\Big)\,.
\eeao
However, this alternative expression of the Laplace \fct al $\Psi_N$ 
corresponds to another \rep\ of the \pp\ $N$:
\beam\label{eq:march11c}
N=\sum_{i=1}^\infty\sum_{j=-\infty}^\infty\vep_{(\Gamma_i/\theta_{|X|})^{ -1/\a} \wt Q_{ij}}\,,
\eeam
where the Poisson points $(\Gamma_i^{-1/\a})$ are independent of the 
\seq\ $\big(\sum_{j\in\bbz}\vep_{\wt Q_{ij}}\big)$ of iid copies of
$\sum_{j\in\bbz} \vep_{\wt Q_j}$. 
\par
We observe that $|\wt Q_j|\le 1$ a.s. and
$|\wt Q_{T^\ast}|=1$ a.s. 
The extremal index $\theta_{|X|}$ plays an important 
role in \rep\ \eqref{eq:march11c}. Each Poisson point 
$(\Gamma_i/\theta_{|X|})^{ -1/\a}$ stands for the radius of a circle
around the origin, and the points $(\wt Q_{ij})_{j\in\bbz}$ are inside or
on this circle. In this sense, each Poisson point 
$(\Gamma_i/\theta_{|X|})^{ -1/\a}$ creates an  extremal cluster. Therefore
we refer to the process $N$ as a {\em cluster Poisson process.}

\subsubsection{Equivalent expressions for the extremal index}
Based on the results in the previous subsection we can derive
equivalent expressions of $\theta_{|X|}$ in terms of $Q_{T^\ast}$ and $T^\ast$.
\bpr\label{cor:cluster1}
Assume the conditions of Theorem~\ref{thm:davishsing2}.
Then the extremal index $\theta_{|X|}$ of $(|X_t|)$ 
coincides with the following quantities:
\beam\label{eq:august28b}
%\E\Big[ \sup_{j\in\bbz} |Q_j|^\a \Big]
%=\E\Big[ \dfrac{\sup_{j\in\bbz} |\Theta_j|^\a} {\sum_{t\in\bbz} |\Theta_t|^\a}\Big]\label{eq:august21xx}&=&
%\E\Big[  \dfrac{ |\Theta_{T ^\ast}|^\a}{\|Theta\|_\a^\alpha}\Big]
 \E[|Q_{T^\ast}|^\a]
=\P(Y\,|Q_{T^\ast}|>1)=\P(T^\ast=0)\,.
\eeam
Here $Y$ is a Pareto$(\a)$ \rv\ independent
of $Q_{T^\ast}$ and $T^\ast$ is the time of the largest record of $(|\Theta_t|)$.
\epr
\bre
We observe that 
\beao
 \E[|Q_{T^\ast}|^\a] = \E\Big[\dfrac{\max_{t\in\bbz} |\Theta_t|^\a}{\sum_{j\in\bbz} |\Theta_j|^\a}\Big]=\theta_{|X|}\,.
\eeao
Since $\theta_{|X|}=\P(T^\ast=0)$ the extremal index $\theta_{|X|}$ 
has the intuitive interpretation
as the \pro y that $(|\Theta_t|)$ assumes its largest value at time zero.
\ere
\bexam\label{ex:br}\rm
We consider the \regvary\ solution of an affine \sre\ under the 
conditions and with the notation of Example~\ref{exam:2}. An exception where the extremal index has an explicit solution
is the case $\log A_t= N_t-0.5$ for an iid standard normal \seq\ $(N_t)$. 
Then $\E[A_t]=1$ and the theory mentioned in Example~\ref{exam:2} yields
\regvar\ of $(X_t)$ with index 1. Using the expression $\P(T^\ast=0)$ and applying some random walk theory (such as the 
results in \cite{chang:peres:1997}), one obtains an exact 
expression for $\theta_X$ in terms of the Riemann zeta \fct\ $\zeta$; 
see Example \ref{ex:br2}. A first order
approximation to this formula is given by
\beam\label{eq:may1a} 
\theta_X \approx \dfrac 1{2}\exp\Big(\dfrac{ \zeta(0.5)}{\sqrt {2\pi}} \Big)\approx \dfrac 1{2}\exp(-0.5826)\approx 0.2792.
\eeam
%{\red Note that this expression remains valid for $\theta_{|X|}$ when $B$ takes negative values as the spectral tail process is preserved.}
\eexam

\bexam\label{ex:br2}\rm 
Let $B^{(i)} = (B_t)_{t \in \R}$ be iid standard Brownian motions independent of 
$\Gamma_1<\Gamma_2<\cdots$ which are the points of a unit-rate Poisson process on 
$(0,\infty)$. 
We consider the stationary max-stable {\em Brown-Resnick \cite{brown:resnick:1977} process}
\beao
X_t = \sup_{i \geq 1} \Gamma_i^{-1}\, \ex^{\sqrt{2} \,B_t^{(i)} - |t|} \,, \qquad t \in \R\,.
\eeao
It has unit Fr\'echet marginals $\P(X_t\le x)=\Phi_1(x)=\ex^{-x^{-1}}$, $x>0$.
Any discretization $X^{(\delta)}=(X_{\delta\,t})_{t\in \Z}$ for $\delta>0$ 
is \regvary\ with index 1 and 
spectral tail process 
$\Theta_t^{(\delta)}=\ex^{\sqrt{2} B_{\delta\, t} - \delta |t|}$, $t\in\bbz$.
Direct calculation of $-x\,\log \P(n^{-1}\,\max_{1\le t\le n} X_{\delta\,t}\le x)$, $x>0$, yields the extremal index of $X^{(\delta)}$ as the limit
\beqq\label{eq:extr_index_BM}
\theta_X^{(\delta)}
=    \lim_{n \to \infty} n^{-1} \E\Big[ \sup_{0\le t\le n} \ex^{\sqrt{2} B_{\delta\,t} - \delta t}\Big]\,.
\eeqq
We use the expression $\theta_X^{(\delta)}=\P(T^{\ast(\delta)}=0)$ where 
$T^{\ast(\delta)}$ is the first record time of $(\Theta_t^{(\delta)})_{t\in\bbz}$; see \eqref{eq:august28b}.
We consider the first ladder height epoch
$\tau_+(\delta)=\inf\{t\ge 1\,: \sqrt 2\, B_{\delta\, t} +\delta t<0\}$. 
Using the symmetry of the Gaussian \ds ,    $(\Theta_t^{(\delta)})_{t\ge 1}\eqd
 (1/\Theta_{-t}^{(\delta)})_{t\ge 1}$, we obtain 
$\theta_X^{(\delta)}=\P(T^{\ast(\delta)}=0)=\P(\tau_+(\delta)=\infty)^2$. Combining this with the classical identity $\P(\tau_+(\delta)=\infty)=1/{\E[\tau_-(\delta)]}$ for $\tau_-(\delta)=\inf\{t\ge 1\,:\sqrt 2\, B_{\delta\, t} -\delta t\le 0\}$, from random walk theory (see \cite{asmussen:2003}) we get
\beao 
\theta^{(\delta)}_X = \Big(\dfrac{1}{\E[\tau_-(\delta)]}\Big)^2 = \Big(\dfrac{\E[B_\delta -\delta]}{\E[\sqrt 2 B_{\tau_-(\delta)}-\tau_-(\delta)]}\Big)^2=\delta^2(\E[\sqrt{2}\,B_{\tau_+(\delta)}+\tau_+(\delta)])^{-2}\,,
\eeao
where we used Wald's lemma and the symmetry of the Gaussian distribution. To be able to apply Theorem 1.1 in
\cite{chang:peres:1997} we standardize the increments of the random walk $\sqrt 2 B_{\delta\,t}$ dividing them by $\sqrt{2\delta}$, turning the drift into $\sqrt {\delta/2}$, and we get
\beao 
\E[\sqrt 2 \,B_{\tau_+(\delta)}+\tau_{+}(\delta)]= \sqrt{ \delta}  \exp\Big (-\dfrac{\sqrt \delta}{2\sqrt \pi} \sum_{n=0}^\infty \dfrac{\zeta(1/2-n)}{n!(2n+1)} \Big(-\dfrac{\delta}4 \Big)^n\Big)\,.
\eeao
This implies that
\beao
\theta_X^{(\delta)} =\delta \exp\Big (\sqrt{\dfrac{  \delta}{ \pi}} \sum_{n=0}^\infty \dfrac{\zeta(1/2-n)}{n!(2n+1)} \Big(-\dfrac{\delta}4 \Big)^n\Big).
\eeao
We recover the {\em Pickands constant} of the 
Brown-Resnick process (see \cite{pickands:1969}) as the limit $\lim_{\delta\downarrow 0}\delta^{-1}\, \theta_X^{(\delta)}$:
\beao
\mathcal H_X^{(0)}=\lim_{T\to \infty}\dfrac1T\E\Big[ \sup_{0\le t\le T} \ex^{\sqrt{2} B_{t} - t}\Big]=1.
\eeao 

\eexam

\begin{proof}[Proof of Proposition \ref{cor:cluster1}]
Consider the supremum of all points of the limit process $N$
in Theorem~\ref{thm:davishsing2}:
\beao
M=
\sup_{i\ge 1}\,\Gamma_i^{-1/\a} \sup_{j\in\bbz} |Q_{ij}|\,.
\eeao
The \seq s $(\Gamma_i)$ and $(Q^{(i)})$  are independent and 
$M=\sup_{i\ge 1}\,\Gamma_i^{-1/\a} V_i$
for the iid \seq\ $V_i:= \sup_{j\in\bbz} |Q_{ij}|$, $i=1,2,\ldots,$ 
whose generic element $V$ has the
property $\E[V^\a]<\infty$. Indeed, $V\le 1$ a.s. by construction.
The points $(\Gamma_i^{-1/\a},V_i)$ constitute a marked Poisson process $N_{\Gamma,V}$
with state space $E=(0, \infty)\times [0,\infty)$ and mean \ms\ given by 
$\mu((x,\infty)\times [0,y])= x^{-\a}\,F_V(y)$, $x>0,y\ge 0$, where 
$F_V$ is the \df\ of $V$.
For $x>0$ we consider 
$B_x= \{(y,v)\in E: y\,v>x\}$. We observe that
\beao
\mu(B_x)= \int_{v=0}^\infty\int_{y=x/v}^\infty\a \,y^{-\a-1}\,F_V(dv)=
\int_0^\infty (x/v)^{-\a}\,F_V(dv)= x^{-\a}\,\E[V^\a]\,.
\eeao
Therefore 
we have for $x>0$,
\beao
\P(M\le x)&=&\P\big(\Gamma_i^{-1/\a} \,V_i\le x\,,i\ge 1\big)\\
&=&P(N_{\Gamma,V}(B_x)=0)\\
&=&\ex^{-\mu(B_x)}= \ex^{-x^{-\a}\,\E[V^\a]}\,.
\eeao
Thus $M$ is a scaled version of the  standard Fr\'echet \ds , 
$\Phi_\a(x)=\ex^{-x^{-\a}}$, $x>0$:
\beao
 \P(M\le x)
&=&\Phi_\a^{\E[V^\alpha]}(x)\,,\qquad x>0\,.
\eeao 
On the other hand,   Theorem~\ref{thm:davishsing2} and an 
application of the \cmt\ yield as $\nto$,
\beao
\P\big(a_n^{-1}M_n\le x\big)&=&\P\big(N_n(x,\infty)=0\big)\\
&\to &\P\big(N(x,\infty)=0\big)\\
&=&\P(M\le x)\,,\qquad x>0\,.
\eeao
In view of the definition of the extremal index of the \seq\ $(|X_t|)$
we can identify 
\beao%\label{eq:august21xx}
\E[V^\a]&=&\E\Big[ \sup_{j\in\bbz} |Q_j|^\a \Big]=\E[|Q_{T^\ast}|^\a].
%= \E\Big[ \dfrac{\sup_{j\in\bbz} |\Theta_j|^\a} {\sum_{t\in\bbz} |\Theta_t|^\a}\Big]
\eeao
 as the  value $\theta_{|X|}$. 
This proves the first part of 
\eqref{eq:august28b}.
The identity
\beao
\E[|Q_{T^\ast}|^\a]=
\P(Y\,|Q_{T^\ast}|>1)=\P(|Q_{T^\ast}|^\a>Y^{-\a}).
\eeao
is immediate since $Q$ and $Y$ are independent, and $Y^{-\a}$ is $U(0,1)$ distributed. 
\par
Applying the time-change formula \eqref{eq:tcp}, shifting $k$ to zero, we obtain
\beao
\theta_{|X|}&=&\E[|Q_{T^\ast}|^\a]\\
&=&\sum_{k\in \Z}\E\Big[\dfrac{ |\Theta_{k}|^\a}{\sum_{j\in\bbz}|\Theta_j|^\alpha}\, \1(T ^\ast=k)\Big]\\
&=&\sum_{k\in \Z}\E\Big[  \dfrac{ |\Theta_{-k}|^\a}{\sum_{j\in\bbz}|\Theta_{j-k}|^\alpha}\,  \1(T ^\ast=0)\Big]\\
&=&\P(T ^\ast=0)\,.
\eeao
This proves  the last identity in \eqref{eq:august28b}.
\end{proof}

\section{Estimation of  the extremal index - a short review and a new estimator based 
on the spectral cluster process}
First approaches to the estimation of the extremal index are due to 
\cite{hsing:1993,smith:weissman:1994}. 
Estimators based on exceedences of a threshold were proposed in 
\cite{ferro:segers:2003,suveges:2007,suveges:davison:2010,ferro:segers:robert:2009}.  A modern approach to the maxima method was started  
in \cite{northrop:2015}; improvements and asymptotic limit theory can be found in \cite{berghaus:bucher:2018,bucher:jennessen:2020}.
\par 
We will consider some  standard estimators of $\theta_{X}$.
For the sake of argument  we assume that $(X_t)$ is a non-negative 
stationary  process with marginal
\ds\ $F$, $k_n=n/r_n$ is an integer \seq\ \st\ $r_n\to\infty$,
$k_n\to\infty$, and $(u_n)$ is a threshold \seq\ satisfying $u_n\uparrow x_F$.

\subsection{Blocks estimator}
Recall that $\theta_X$ has interpretation as the reciprocal of the 
expected size of extremal clusters. This idea is the basis for 
inference procedures from the early 1990s (see 
\cite{smith:1989,davison:smith:1990}). 
Clusters are identified as blocks of length $r=r_n$ with at least one exceedance of a high threshold $u=u_n$. A blocks estimator $\wh\theta^{\rm bl}$ is given by 
the ratio of the number $K_n(u)$ of such clusters and the total number of exceedences $N_n(u)$:
\begin{align}\label{eq:cluster:estimator}
    \widehat{\theta}^{\rm bl}_u(r) = \frac{K_n(u)}{N_n(u)} := \frac{ \sum_{t=1}^{k_n} \1( M_{(t-1)r +1,t\,r} > u ) }{ \sum_{t=1}^n \1(X_t > u) }.
\end{align}
This method requires the choice of block length $r$ 
and  threshold level $u $  satisfying $r_n\overline{F}(u_n) \to 0$; 
if $r_n\to\infty$ does not hold at the  prescribed  rate $\wh\theta^{\rm bl}$
is biased. Estimators using clusters of extreme exceedences were also
considered in \cite{hsing:1993}.
\par
A slight modification of the blocks estimator is the 
{\em disjoint blocks estimator} of \cite{smith:weissman:1994}: 
\beao
 \widehat{\theta}^{\rm dbl} = \dfrac{ \log(1-K_n(u)/k_n)}{r\, \log(1-N_n(u)/n)}.
\eeao
Assuming some weak dependence condition on $(X_t)$, the heuristic  idea behind the estimator is the approximations
\beao
\big(\P(M_r\le u_n)\big)^{k_n}\approx \P(M_n\le u_n)\approx F^{\theta_X\,n}(u_n),
\eeao 
for a suitable \seq\ $(u_n)$. Then, taking logarithms and replacing
$\ov F(u_n)$ and $\P(M_n>u_n)$ by their empirical estimators $N_n(u)/n$ and 
$K_n(u)/k_n$, respectively, we obtain
\beao
\theta_X&\approx& \dfrac{\log \P(M_n\le u_n)}{n\,\log F(u_n)}=
\dfrac{\log (1-\P(M_n>u_n))}{n\,\log(1-\ov F(u_n))}\\
&\approx &\dfrac{\log (1-K_n(u)/k_n)}{r_n\,\log (1- N_n(u)/n)}=\wh\theta^{\rm dbl}\,.
\eeao
Assuming that both $K_n(u)/k_n$ and $N_n(u)/n$ converge to zero,
a Taylor expansion of $\log (1+x)=x(1+o(1))$ as $x\to 0$ shows that
$\wh \theta^{\rm bl}\approx \wh\theta^{\rm dbl}$.  \cite{smith:weissman:1994} showed that   $\wh\theta^{\rm dbl}$ has a smaller \asy\ variance than $\wh \theta^{\rm bl}$. \cite{ferro:segers:robert:2009} proposed a sliding blocks version of $\wh \theta^{\rm dbl}$ with an even smaller \asy\ variance.
\beam\label{eq:sliding}
    \widehat{\theta}^{\rm slbl}(u,r) =  \frac{ - \log\Big( \frac{1}{n-r+1} \sum_{t=1}^{n-r+1}\1(M_{t,t+r} \le u ) \Big)}{ N_n(u)/k_n }\,.
\eeam

 \subsection{Runs and intervals estimator}
\cite{smith:weissman:1994} proposed the alternative 
\textit{runs estimator}. It is based on the limit relation
\eqref{eq:obrien}: the \pro y $\P(M_{\ell_n}\le u_n\,\mid \,X_0>u_n)$ is replaced by a sample version for some \seq\ 
$l=l_n\to\infty $:
\beam\label{eq:runs}
\wh \theta^{\rm runs}_u(l)=\frac{1}{N_n(u)}\sum_{i=1}^{n-l}\1(X_i>u_n\,,M_{i+1,i+l}\le u_n)\,.
\eeam
Clusters are considered distinct if 
they are separated by at least $l$ observations not exceeding $u$.
In \cite{ferro:segers:2003} a complete study of the runs estimator and the 
inter-exceedence times is given. The thresholds $(u_n)$ need to satisfy
$r_n\ov F(u_n)\to 1$, and $l_n\le r_n$.
\par
Consider the {\em exceedance times}:
\beao
S_0(u)=0\,, \qquad S_i(u) = \min \{ t > S_{i-1}(u) : X_t > u_n\}\,, \qquad i\ge 1\,,
\eeao
with {\em inter-exceedance times} $T_i(u)=S_i(u)-S_{i-1}(u)$, 
$i\ge 1$. The \seq\ $(T_i(u))_{i\ge 2}$ constitutes a stationary \seq .
If $r_n\,\overline{F}(u_n) \to 1$, \cite{ferro:segers:2003} noticed that  
$(n\,T_2(u))$ converges in \ds\ to a limiting mixture given by 
$(1-\theta_{X})\1_{0}(x)+
\theta_X\,(1-\ex^{-\theta_X\,x})$, $x\ge 0$. Calculation yields to the 
coefficient of variation $\nu$ of $T_2(u)$ whose square is given by 
\beao
\nu^2=\var(T_2(u))/(\E[T_2(u)])^2=\E[T_2^2(u)]/(\E[T_2(u)])^2-1= 2/\theta_X-1\,,
\eeao
leading to overdispersion $\nu>0$
\fif\ $\theta_X<1$. Replacing the moments on the \lhs\ by sample versions and adjusting the empirical moments for bias,
\cite{ferro:segers:2003} arrived at the {\em intervals estimator} 
\beam\label{eq:intervalts}
\widehat{\theta}^{\rm int}(u) &=& 
%    1\wedge  \dfrac{2 \big(\sum_{i=2}^{N_n(u)} T_i(u)\big)^2 }{(N_n(u)-1) \sum_{i=2}^{N_n-1} T_i^2(u)} &\text{if} \max\{T_i\} \le 2\\
     1\wedge   \tfrac{2 \big(\sum_{i=2}^{N_n(u)} (T_i(u)-1)\big)^2 }{(N_n(u)-1) \sum_{i=2}^{N_n(u)} (T_i(u)-1)(T_i(u) - 2)}\,.
\eeam
See also \cite{suveges:2007,suveges:davison:2010}.

\subsection{Northrop's estimator}
Assume for the moment that $(X_i)$ is iid and $F$ is continuous. Then $F(X)$ is uniform on $(0,1)$.
Hence for $r=r_n$ and $x>0$,
\beao
\P\big(-r_n\,\log F(M_r)>x\big)&=& \P(F(M_r)\le  \ex^{-x/r})\\&=& 
\P(\max_{i=1,\ldots,r_n} F(X_i)\le  \ex^{-x/r})\\&=&\big(\P(F(X)\le \ex^{-x/r})\big)^r=\ex^{-x}\,.
\eeao
For a weakly dependent \seq\ $(X_i)$ with marginal \ds\ $F$, assume the existence of the extremal 
index for $(F(X_t))$ which, by monotonicity of $F$, coincides with $\theta_X$:
\beao
\P\big(-r_n\,\log F(M_r)>x\big)=\P(\max_{i=1,\ldots,r_n} F(X_i)< \ex^{-x/r})
\to \ex^{-\theta_X\,x}\,, \qquad x>0\,.
\eeao
Thus the \rv s $(-r_n\,\log F(M_r))$ are \asy ally ${\rm Exp}(\theta_X)$ 
distributed. For iid ${\rm Exp}(\theta_X)$ \rv s the maximum likelihood 
estimator of $\theta_X$ is given by the reciprocal of the sample mean.
These ideas lead to Northrop's estimators 
\cite{northrop:2015}. Mimicking the maximum likelihood estimator of iid ${\rm Exp}(\theta_X)$  data 
for a stationary \seq\ $(X_t)$, one considers the quantities $-r_n\,\log F(M_{t,t+r})$, $t=1,\ldots,n-r_n$, and constructs sliding 
or disjoint blocks estimators of $\theta_X$:
\beam\label{eq:Northrop}
    \widehat{\theta}^{\rm Nsl}(r) &=& \Big( \frac{1}{n-r+1}\sum_{t=1}^{n-r + 1} (- r\,  \log F_n( M_{t,t+r} ))\Big)^{-1}\,,\\
 \widehat{\theta}^{\rm Ndbl}(r) &=& \Big( \frac{1}{[n/r]}\sum_{i=1}^{[n/r]} (- r\,  \log F_n( M_{r\,(i-1)+1,r\,i} ))\Big)^{-1}\,.
\eeam
Here $F_n$ is the empirical \df\ of the data. This particular choice 
of estimator of $F$ depends on the whole sample, hence introduces additional dependence. This fact requires an optimal choice of block length $r_n$ for implementation.

\subsection{An estimator based on the spectral cluster process}
In this subsection we consider a stationary non-negative 
\regvary\ process $(X_t)$ with index $\a>0$,
spectral tail process $(\Theta_t)$ and normalizing \seq\ $(a_n)$ satisfying 
$n\,\P(X>a_n)\to 1$.   Proposition~\ref{cor:cluster1} yields the alternative \rep\  
$\theta_X=\E[Q_{T^\ast}^{\alpha}] $ where $(Q_t)$ is the spectral cluster process of $(X_t)$. We will construct an estimator based on this identity.
\par
We consider sums and maxima over disjoint blocks of size $r=r_n=o(n)$:
\beao
S_{i,r}^{(\alpha)} := \sum_{t=(i-1)\,r + 1}^{i \,r} X_t^\alpha\,,\qquad M_{i,r}=\max_{t=(i-1)r+1,\ldots,i\,r} X_t\,, \qquad i=1,\ldots,k_n\,.
\eeao
The following limit relation is proved in 
\cite{buritica:mikosch:wintenberger:2021}:
\beam\label{eq:limit:theta}
\lim_{n \to  \infty} \E\big[ \,  M_{1,r}^\alpha/S_{1,r}^{(\alpha)} \, | \, S_{r}^{(\alpha)} > a_n^\alpha \big] = \E[Q_{T^\ast}^\alpha],
\eeam
which is based on \ld\ results for \regvary\ stationary \seq s; see for example \cite{buritica:mikosch:wintenberger:2021}.
Now we build an estimator of $\theta_X$ from an empirical version of the left-hand expectation.
Define the corresponding estimator by 
 \begin{align}\label{eq:estimate:thet2}
     \widehat{\theta}_v^{\rm scp}(r) := \frac{  \sum_{i=1}^{k_n}  \frac{M_{i,r}^\alpha}{S^{(\alpha)}_{i,r}}\,  \1\big( S_{i,r}^{(\alpha)}> v \big)}{\sum_{i=1}^{k_n}    \1\big( S_{i,r}^{(\alpha)}> v \big)}\,.
\end{align}
Here we choose $v=S_{(s),r}^{(\alpha)}$,
the $s$th largest among $(S_{i,r}^{(\alpha)})_{i=1,\ldots,k_n}$ for an integer \seq\ 
$s=s_n$ \st\  $s_n=o(k_n)$.

\section{A Monte-Carlo study of the estimators}
We run a short study based on $1\,000$ simulated
processes   $(X_t)_{t=1,\ldots,5000}$ 
for comparing the performances of 
some of the aforementioned estimators. First, $(X_t)$ is an AR(1) process
with parameter $\varphi=0.2$ and iid student(1) noise,
resulting in a \regvary\ process with index 1 and $\theta_{|X|}=0.8$. 
Second, we consider the 
\regvary\ solution of an affine  \sre\ with iid 
$\log A_t\sim N(-0.5,1)$, $B_t\equiv 1$, %N(0,1)$ 
and $\theta_{X}\approx  0.2792$; see \eqref{eq:may1a}.
\par
\begin{figure}[!htb]
    \centering
    
   \includegraphics[height=11.5cm,width=11.5cm]{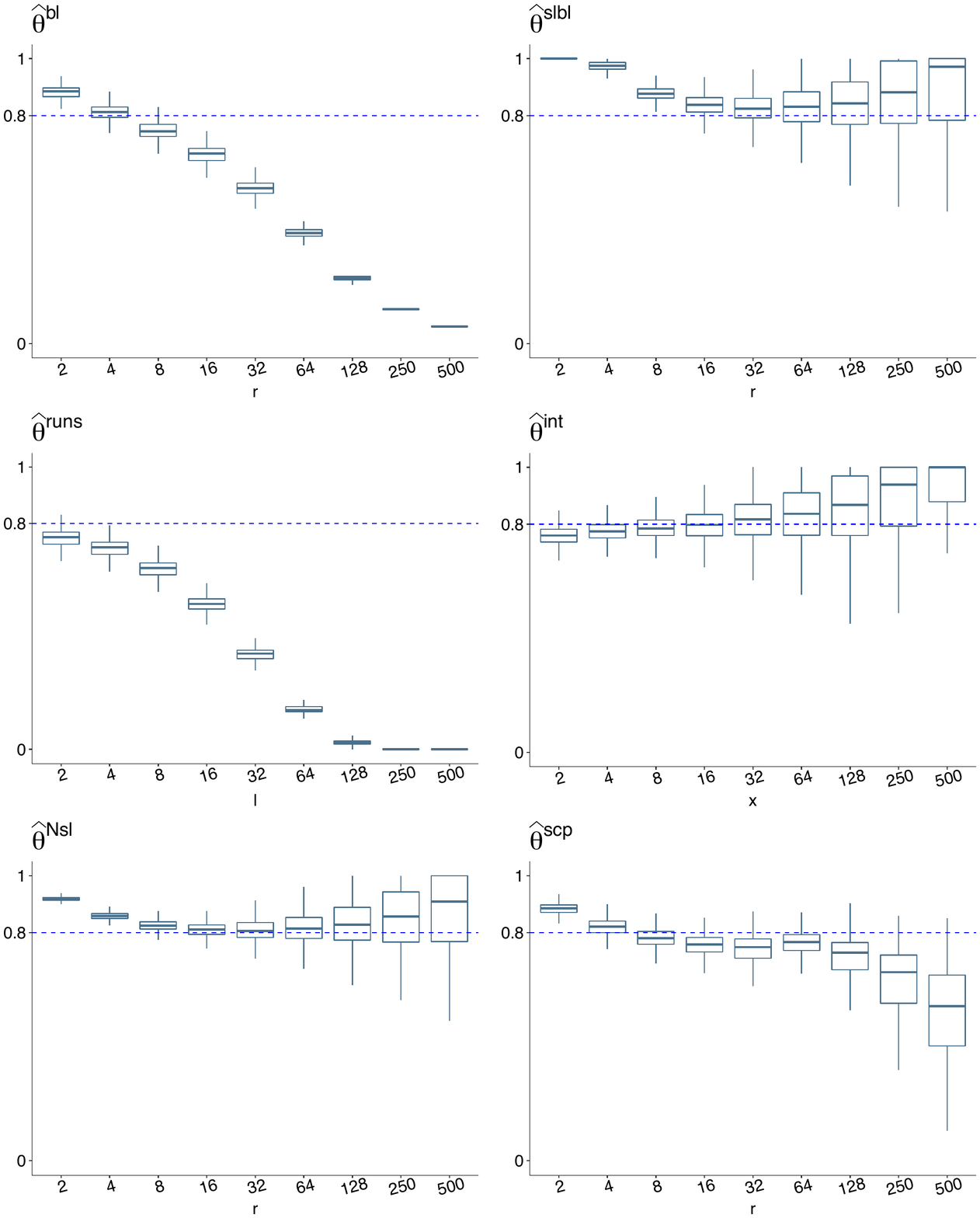}
    
    \bfi{\small  Boxplots based on $1\,000$ simulations for the estimation of 
    $\theta_{|X|}=0.8$ in the $AR(1)$ model with $\varphi=0.2$ and iid student$(1)$ noise.    }\label{fig:AR}\efi
\end{figure}
\begin{figure}[!htbp]
\centering
   \includegraphics[width=11.5cm,height=11.5cm]{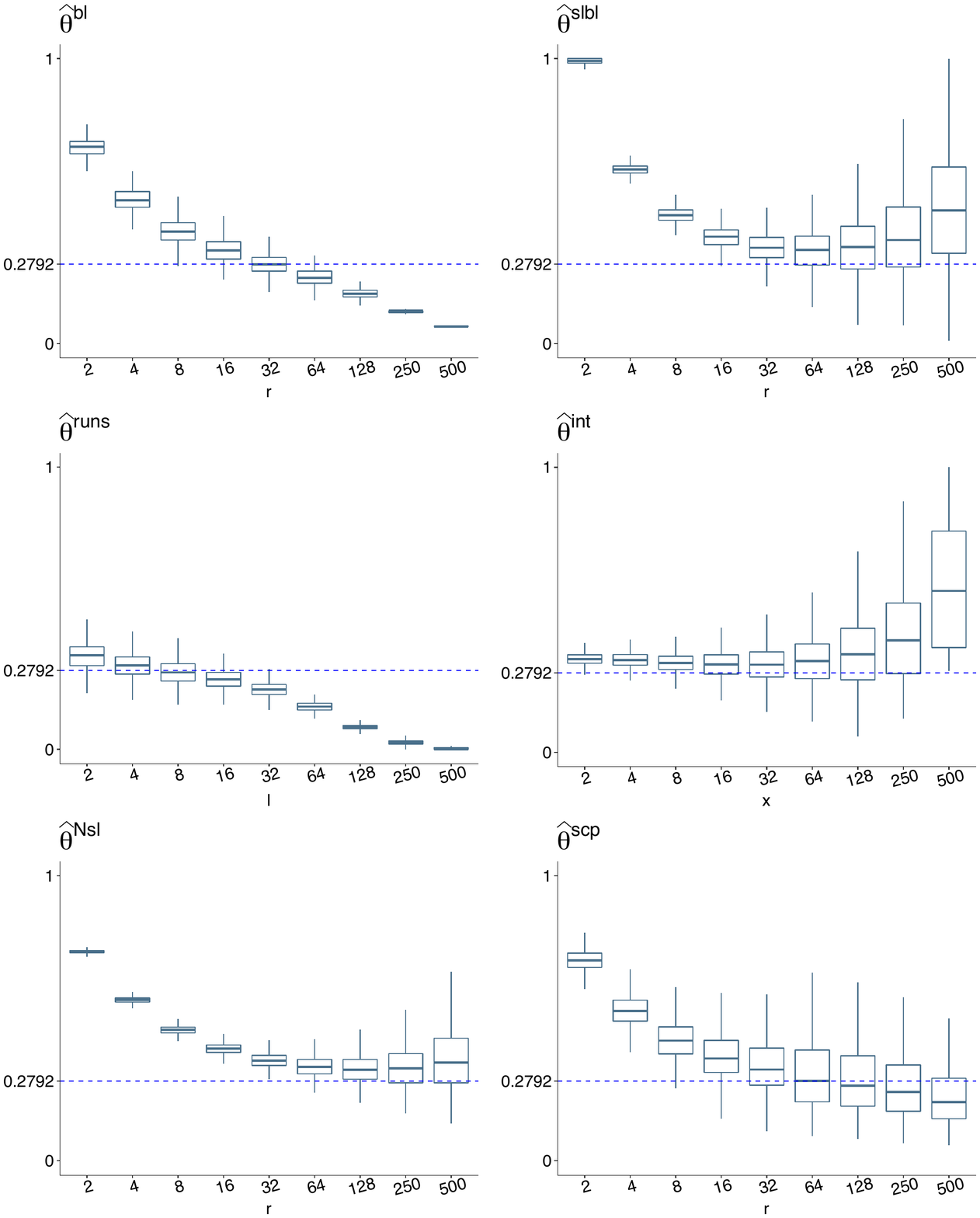}
\bfi{\small   Boxplots based on $1\,000$ simulations for the estimation of  $\theta_{X}\approx 0.2792$ for the solution to a \sre .
}\label{fig:ARCH}\efi 
\end{figure}

Figures~\ref{fig:AR} and \ref{fig:ARCH} show boxplots of the simulation study. 
\begin{itemize}
\item
$\widehat{\theta}^{\rm bl}$ and $\widehat{\theta}^{\rm runs}$ 
are \fct s of the block and run lengths, 
respectively. $u$ is the largest  $[ n^{0.6} ]$th 
upper order statistic of the sample. 
\item
$\widehat{\theta}^{\rm slbl}$ is a function of $r$. $u$ is the $r$th upper order statistic of the sample.
\item
$\widehat{\theta}^{\rm int}$ is a function  of $x$.  
$u$ is the $[n/x]$th upper order statistic.
\item
$\widehat{\theta}^{\rm Nsl}$, $\widehat{\theta}^{\rm scp}$ are \fct s of $r$.
\item
For $\widehat{\theta}^{\rm scp}$ we choose $s=[n^{0.6}/r]$. The tail index $\a$ 
is estimated by the Hill estimator from  \cite{dehaan:mercadier:zhou:2016} based on $[n^{0.8}]$ upper order statistics of the sample.
\end{itemize}
According to the folklore in the literature, Northrop's estimator $\widehat{\theta}^{\rm Nsl}$ outperforms the classical estimators (runs, blocks); it has smallest variance 
but it may be difficult to control its bias. Our experience with 
$\widehat{\theta}^{\rm scp}$ shows that it performs better 
than the other estimators as regards the bias, especially when $\theta_{X}$ 
is small. 
The intervals estimator $\widehat{\theta}^{\rm int}$ is preferred 
by practitioners because the  choice of the hyperparameter $x$ is robust \wrt\ 
different values of $\theta_{X}$. This cannot be said about the other estimators with the exception of $\widehat{\theta}^{\rm scp}$. 
In our experiments with sample size 
$n=5000$, the choices $x=32$ and $r=64$ work well for   $\widehat{\theta}^{\rm int}$ and $\widehat{\theta}^{\rm scp}$, respectively. We did not fine-tune
the hyperparameter $s$ in $\widehat{\theta}^{\rm scp}$ in our experiments. 

\par
 
\newpage
%\section{Bibliography}

\end{document}